\documentclass[11pt, reqno]{amsart}

\usepackage[OT2, T1]{fontenc}

\usepackage[usenames,dvipsnames]{color}

\newcommand{\revise}[1]{{\color{Bittersweet}   [#1]}}
\newcommand{\ready}[1]{{\color{Brown}   [#1]}}

\renewcommand{\ready}[1]{#1}	    

\usepackage{url}
\usepackage{amsmath}
\usepackage{array}
\usepackage{graphicx}
\usepackage{amssymb}
\usepackage{amsthm}
\definecolor{link}{RGB}{11,0,128}
\usepackage[colorlinks=true,citecolor=NavyBlue,linkcolor=Bittersweet,urlcolor=link]{hyperref}
\usepackage{hyperref}
\usepackage{paralist}
\usepackage{colonequals}		
\usepackage{color}
\usepackage{mathabx}		
\usepackage{amscd}
\usepackage[all,cmtip]{xy}	
\usepackage{sseq}			
\usepackage{verbatim}
\usepackage{parskip}
\usepackage{microtype}
\usepackage[in]{fullpage}
\usepackage{mathrsfs} 			
\usepackage{cleveref}
\usepackage{enumitem}
\usepackage{moreenum}			
\usepackage[alphabetic,lite]{amsrefs} 	
\usepackage{stmaryrd}	
\usepackage[foot]{amsaddr}
\usepackage{subfiles}
\usepackage{ragged2e}

\DeclareSymbolFont{cyrletters}{OT2}{wncyr}{m}{n}
\DeclareMathSymbol{\Sha}{\mathalpha}{cyrletters}{"58}

\newcommand{\gA}{\alpha}

\newcommand{\bA}{\mathbb{A}}
\newcommand{\bC}{\mathbb{C}}

\newcommand{\bF}{\mathbb{F}}
\newcommand{\bG}{\mathbb{G}}

\newcommand{\bP}{\mathbb{P}}
\newcommand{\bQ}{\mathbb{Q}}
\newcommand{\bR}{\mathbb{R}}

\newcommand{\bZ}{\mathbb{Z}}

\newcommand{\cH}{\mathcal{H}}

\newcommand{\cO}{\mathcal{O}}

\newcommand{\fm}{\mathfrak{m}}

\newcommand{\sO}{\mathscr{O}}

\newcommand{\ra}{\rightarrow}

\newcommand{\Ra}{\Rightarrow}

\newcommand{\xra}{\xrightarrow}

\newcommand{\hra}{\hookrightarrow}

\newcommand{\wh}{\widehat}

\newcommand{\ce}{\colonequals}
\newcommand{\ov}{\overline}

\renewcommand{\b}{\textbf}
\newcommand{\surjects}{\twoheadrightarrow}

\newcommand{\tensor}{\otimes} 		
\newcommand{\isomto}{\overset{\sim}{\longrightarrow}}

\newcommand{\gp}{{\mathrm{gp}}}		
\newcommand{\fppf}{\mathrm{fppf}}		
\newcommand{\et}{\mathrm{\acute{e}t}}	
\newcommand{\llb}{\llbracket}		
\newcommand{\rrb}{\rrbracket}		
\newcommand{\lb}{[}		
\newcommand{\rb}{]}		
\newcommand{\tors}{\mathrm{tors}}		
\newcommand{\sh}{\mathrm{sh}}		

\providecommand{\p}[1]{\left(#1\right)}
\providecommand{\up}[1]{{\upshape(}#1{\upshape)}}
\providecommand{\uref}[1]{{\upshape\ref{#1}}}

\providecommand{\ucolon}{{\upshape:} }
\providecommand{\uscolon}{{\upshape;} }

\providecommand{\f}[2]{\frac{#1}{#2}}
\DeclareMathOperator{\Ker}{Ker}			
\DeclareMathOperator{\Spec}{Spec}		
\DeclareMathOperator{\Spa}{Spa}		
\DeclareMathOperator{\Hom}{Hom}			
\DeclareMathOperator{\Char}{char}		
\DeclareMathOperator{\PGL}{PGL}			
\DeclareMathOperator{\Gal}{Gal}	
\DeclareMathOperator{\Norm}{Norm}		
\DeclareMathOperator{\Res}{Res}		
\DeclareMathOperator{\GL}{GL}		
\DeclareMathOperator{\Iso}{Iso}		
\DeclareMathOperator{\rk}{rk}		
\DeclareMathOperator{\Pic}{Pic}		
\newcommand{\ba}{\begin{aligned}}
\newcommand{\ea}{\end{aligned}}
\newcommand{\be}{\begin{equation}}
\newcommand{\ee}{\end{equation}}
\newcommand{\pf}{\begin{proof}}
\newcommand{\bpf}{\begin{proof}}
\newcommand{\epf}{\end{proof}}
\newcommand{\bthm}{\begin{thm}}
\newcommand{\ethm}{\end{thm}}

\newcommand{\bprop}{\begin{prop}}
\newcommand{\eprop}{\end{prop}}
\newcommand{\bcor}{\begin{cor}}
\newcommand{\ecor}{\end{cor}}

\newcommand{\brem}{\begin{rem}}
\newcommand{\erem}{\end{rem}}

\newcommand{\brems}{\begin{rems} \hfill \begin{enumerate}[label=\b{\thenumberingbase.},ref=\thenumberingbase]}
\newcommand{\remi}{\addtocounter{numberingbase}{1} \item}
\newcommand{\erems}{\end{enumerate} \end{rems}}

\newcommand{\begs}{\begin{egs} \hfill \begin{enumerate}[label=\b{\thenumberingbase.},ref=\thenumberingbase]}

\newcommand{\eegs}{\end{enumerate} \end{egs}}

\newcommand{\eremstweak}{\end{enumerate} \end{rems-tweak}}
\newcommand{\eremst}{\end{enumerate} \end{rems-tweak}}
\newcommand{\blem}{\begin{lemma}}
\newcommand{\elem}{\end{lemma}}

\newcommand{\bconj}{\begin{conj}}
\newcommand{\econj}{\end{conj}}
\newcommand{\bprob}{\begin{Problem}}
\newcommand{\eprob}{\end{Problem}}

\newcommand{\bq}{\begin{Q}}
\newcommand{\eq}{\end{Q}}
\newcommand{\benum}{\begin{enumerate}[label={{\upshape(\alph*)}}]}
\newcommand{\benuma}{\begin{enumerate}[label={{\upshape(\arabic*)}}]}
\newcommand{\benumr}{\begin{enumerate}[label={{\upshape(\roman*)}}]}
\newcommand{\eenum}{\end{enumerate}}
\newcommand{\bc}{}
\newcommand{\bd}{\begin{defn}}
\newcommand{\ed}{\end{defn}}

\newcommand{\beg}{\begin{eg}}
\newcommand{\eeg}{\end{eg}}

\newcommand{\bcl}{\begin{claim}}
\newcommand{\ecl}{\end{claim}}
\newcommand{\lab}{\label}

\newcommand{\x}{\text}

\newcommand{\q}{\quad}
\newcommand{\qq}{\quad\quad}
\newcommand{\qqq}{\quad\quad\quad}

\newcommand{\tst}{\textstyle}

\newcommand*{\QED}{\hfill\ensuremath{\qed}}

\usepackage{aliascnt}
\newaliascnt{numberingbase}{subsection}

\theoremstyle{plain}
\newtheorem{thm}[numberingbase]{Theorem}
\Crefname{thm}{Theorem}{Theorems}

\Crefname{rethm}{Theorem}{Theorem}
\newtheorem{prop}[numberingbase]{Proposition}
\Crefname{prop}{Proposition}{Propositions} 
\newtheorem{Q}[numberingbase]{Question}
\Crefname{Q}{Question}{Questions}
\newtheorem{Problem}[subsection]{Problem}
\Crefname{Problem}{Problem}{Problems}
\newtheorem{conj}[numberingbase]{Conjecture}
\Crefname{conj}{Conjecture}{Conjectures}
\newtheorem{cor}[numberingbase]{Corollary}
\Crefname{cor}{Corollary}{Corollaries}
\newtheorem{lemma}[numberingbase]{Lemma}

\Crefname{subprop}{Proposition}{Propositions}

\Crefname{subcor}{Corollary}{Corollaries}

\Crefname{sublem}{Lemma}{Lemmas}

\theoremstyle{remark}
\newtheorem{claim}[equation]{Claim}
\Crefname{claim}{Claim}{Claims}

\Crefname{subrem}{Remark}{Remarks}

\theoremstyle{definition}
\newtheorem{defn}[numberingbase]{Definition}
\Crefname{defn}{Definition}{Definitions}

\Crefname{conv}{Convention}{Conventions}
\newtheorem{eg}[numberingbase]{Example}
\Crefname{eg}{Example}{Examples}
\newtheorem{rem}[numberingbase]{Remark}
\Crefname{rem}{Remark}{Remarks}
\newtheorem*{rems}{Remarks}
\newtheorem*{egs}{Examples}

\theoremstyle{plain}
\newtheorem{thm-tweak}[subsection]{Theorem}
\Crefname{thm-tweak}{Theorem}{Theorems}
\newtheorem{lemma-tweak}[subsection]{Lemma}
\Crefname{lemma-tweak}{Lemma}{Lemmas}
\newtheorem{cor-tweak}[subsection]{Corollary}
\Crefname{cor-tweak}{Corollary}{Corollaries}
\newtheorem{prop-tweak}[subsection]{Proposition}
\Crefname{prop-tweak}{Proposition}{Propositions} 
\newtheorem{conj-tweak}[subsection]{Conjecture}
\Crefname{conj-tweak}{Conjecture}{Conjectures} 

\theoremstyle{definition}
\newtheorem{defn-tweak}[subsection]{Definition}
\Crefname{defn-tweak}{Definition}{Definitions}
\newtheorem{eg-tweak}[subsection]{Example}
\Crefname{eg-tweak}{Example}{Examples}
\newtheorem*{rems-tweak}{Remarks}
\newtheorem{rem-tweak}[subsection]{Remark}
\Crefname{rem-tweak}{Remark}{Remarks}

\newtheoremstyle{subsection-tweak}
   {11pt}
   {3pt}%
   {}
   {}%
   {\bfseries}
   {}%
   {.5em}
   {\thmnumber{\@{#1}{}\@{#2}.}%
    \thmnote{~{\bfseries#3.}}}    
    
\theoremstyle{subsection-tweak}
\newtheorem{pp}[numberingbase]{}
\newcommand{\bpp}{\begin{pp}}
\newcommand{\epp}{\end{pp}}

\theoremstyle{subsection-tweak}
\newtheorem{pp-tweak}[subsection]{}

\numberwithin{equation}{numberingbase}

\makeatletter
\def\@tocline#1#2#3#4#5#6#7{
    \begingroup 
    \@ifempty{#4}{%
    }{%
    }%

    \parindent\z@ \leftskip#3\relax \advance\leftskip\@tempdima\relax
    #5\hskip-\@tempdima
      \ifcase #1
       \or\or \hskip 2em \or \hskip 1em \else \hskip 3em \fi%
      #6\nobreak\relax
    \dotfill\hbox to\@pnumwidth{\@tocpagenum{#7}}\par
    \nobreak
    \endgroup
  }
 \def\l@section{\@tocline{1}{0pt}{1pc}{}{}}

\renewcommand{\tocsection}[3]{%
  \indentlabel{\@ifnotempty{#2}{\makebox[1.3em][l]{%
    \ignorespaces#1 \bfseries{#2}.\hfill}}}\bfseries{#3}
    \vspace{1.5pt}}

\renewcommand{\tocsubsection}[3]{%
  \indentlabel{\@ifnotempty{#2}{\hspace*{-0.5em}\makebox[2.1em][l]{%
    \ignorespaces#1#2.\hfill}}}#3
    \vspace{1.5pt}}

\makeatother

\makeatletter
\newcommand\appendix@section[1]{%
  \refstepcounter{section}%
  \orig@section*{Appendix \@Alph\c@section. #1}%
}
\let\orig@section\section
\g@addto@macro\appendix{\let\section\appendix@section}
\makeatother

\begin{document}

\title{Purity for the Brauer group}

\author{K\k{e}stutis \v{C}esnavi\v{c}ius}
\address{Laboratoire de Math\'{e}matiques d'Orsay, Univ.~Paris-Sud, CNRS, Universit\'{e} Paris-Saclay, 91405 Orsay, France}
\email{kestutis@math.u-psud.fr}
\date{\today}
\subjclass[2010]{\ready{Primary 14F22; Secondary 14F20, 14G22, 16K50.}}
\keywords{\ready{Brauer group, \'{e}tale cohomology, perfectoid ring, punctured spectrum, purity.}}


\begin{abstract} \ready{A purity conjecture due to Grothendieck and Auslander--Goldman predicts that the Brauer group of a regular scheme does not change after removing a closed subscheme of codimension $\ge 2$. The combination of several works of Gabber settles the conjecture except for some cases that concern $p$-torsion Brauer classes in mixed characteristic $(0, p)$. We establish the remaining cases by using the tilting equivalence for perfectoid rings. To reduce to perfectoids, we control the change of the Brauer group of the punctured spectrum of a local ring when passing to a finite flat cover.
} \end{abstract}


\maketitle

\hypersetup{
    linktoc=page,     
}

\renewcommand*\contentsname{}
\q\\
\tableofcontents

\section{The purity conjecture of Grothendieck and Auslander--Goldman} \lab{intro}

\ready{
Grothendieck predicted in \cite{Gro68c}*{\S6} that the cohomological Brauer group of a regular scheme $X$ is insensitive to removing a closed subscheme $Z \subset X$ of codimension $\ge 2$. This purity conjecture is known in many cases (as we discuss in detail below), for instance, for cohomology classes of order invertible on $X$, and its codimension requirement is necessary: the Brauer group of $\bA^2_\bC$ does not agree with that of the complement of the coordinate axes (see \cite{DF84}*{Rem.~3}). In this paper, we finish the remaining cases, that is, we complete the proof of the following theorem.

\bthm[\Cref{purity-conj-pf}] \lab{purity-conj}
For a locally Noetherian scheme $X$ and a closed subscheme $Z \subset X$ such that for every $z \in Z$ the local ring $\cO_{X,\, z}$ of $X$ at $z$ is regular of dimension $\ge 2$, we have
\[
H^2_\et(X, \bG_m) \isomto H^2_\et(X - Z, \bG_m) \qq \text{and} \qq H^3_\et(X, \bG_m) \hra H^3_\et(X - Z, \bG_m).
\]
\ethm

The purity conjecture of Grothendieck builds on an earlier question of Auslander--Goldman pointed out in \cite{AG60}*{7.4}. Due to a result of Gabber \cite{Gab81}*{II, Thm.~1}, that is, due to the agreement of the Brauer group of an affine scheme with its cohomological counterpart, a positive answer to their question amounts to the affine case of the following consequence of \Cref{purity-conj}.

\bthm[\Cref{AG-cor-pf}] \lab{AG-cor}
For a Noetherian, integral, regular scheme $X$ with function field $K$,
\[
\tst H^2_\et(X, \bG_m) = \bigcap_{x \in X\,\text{of height}\, 1} H^2_\et(\cO_{X,\, x}, \bG_m) \qq \text{in} \qq H^2_\et(K, \bG_m).
\]
\ethm

The global \Cref{purity-conj,AG-cor} are known to readily reduce to the following key local purity result.

\bthm[\Cref{win}] \lab{punc-thm}
For a strictly Henselian, regular, local ring $R$ of dimension $\ge 2$, 
\[
H^2_\et(U_R, \bG_m) = 0, \qq \text{where $U_R$ is the punctured spectrum of $R$.}
\]
\ethm

In turn, as we now summarize, many cases of \Cref{punc-thm} are already known.
\benumr
\item \lab{known-1}
The case $\dim R = 2$ follows from the equivalence of categories between vector bundles on $\Spec(R)$ and on $U_R$, see 
\cite{Gro68c}*{6.1 b)}.

\item \lab{known-2}
The case $\dim R = 3$ was settled by Gabber in \cite{Gab81}*{I, Thm.~2}.

\item \lab{known-3}
The vanishing of $H^2_\et(U_R, \bG_m)[p^\infty]$ for the primes $p$ that are invertible in $R$ follows from the absolute purity conjecture whose proof, due to Gabber, is given in \cite{Fuj02} or \cite{ILO14}*{XVI} (special cases also follow from earlier \cite{Gro68c}*{6.1}, \cite{SGA4III}*{XVI, 3.7}, \cite{Tho84}*{3.7}).

\item \lab{known-4}
The vanishing of $H^2_\et(U_R, \bG_m)[p^\infty]$ in the case when $R$ is an $\bF_p$-algebra is given by \cite{Gab93}*{2.5} (and, under further assumptions, also by the earlier \cite{Hoo80}*{Cor.~2}).

\item \lab{known-5}
The case when $R$ is formally smooth over a discrete valuation ring is given by \cite{Gab93}*{2.10}.

\item \lab{known-6}
Gabber announced further cases in an Oberwolfach abstract \cite{Gab04}*{Thm.~5 and Thm.~6} whose proofs have not been published: the case $\dim R \ge 5$ and the case when $R$ is of dimension $4$, of mixed characteristic $(0, p)$, and contains a primitive $p$-th root of unity. 
\eenum

For proving \Cref{punc-thm}, we will use its known cases \ref{known-1}--\ref{known-3} but not \ref{known-4}--\ref{known-6}. 

Our proof has two main steps. The first is to show that the validity of \Cref{punc-thm} for $R$ of dimension $\ge 4$ is insensitive to replacing $R$ by a regular $R'$ that is finite flat over $R$. Such a reduction has also been announced in \cite{Gab04}*{Thm.~4}, but our argument seems simpler and gives a more broadly applicable result. More precisely, we argue in \S\ref{finite-flat} that passage to $R'$ is controlled by the $U_R$-points of a certain homogeneous $R$-space $X$, show that $X$ is affine, and then conclude by deducing that $X(U_R) = X(R)$; the restriction $\dim R \ge 4$ comes from using the vanishing of the Picard group of the punctured spectrum of the local complete intersection $R' \tensor_R R'$ that intervenes in reducing to $X$. In comparison, the argument sketched for \emph{loc.~cit.}~uses deformation theory and a local Lefschetz theorem from \cite{SGA2new}*{X} to eventually obtain passage to $R'$ from the known cases of \Cref{punc-thm}. Since the $p$-primary Brauer group of a perfect $\bF_p$-algebra vanishes, the first step suffices in characteristic $p$.

The second step is to use the tilting equivalence of Scholze introduced in \cite{Sch12} (which, in turn, is a version of the almost purity theorem of Faltings \cite{Fal02}) to show that for a $p$-torsion free perfectoid ring $A$, the $p$-primary cohomological Brauer group $H^2_\et(A[\f{1}{p}], \bG_m)[p^\infty]$ vanishes (see \Cref{p-oid-Br}). This vanishing ultimately comes from the fact that the \'{e}tale $p$-cohomological dimension of an affine, Noetherian scheme of characteristic $p$ is $\le 1$. The intervening comparisons between the \'{e}tale cohomology of (non-Noetherian) affinoid adic spaces and of their underlying coordinate rings add to the technical details required for the second step but not much to the length of the overall argument because, modulo limit arguments, the comparisons we need were proved by Huber in \cite{Hub96}.

The flexibility of the first step leads to the following refinement of \Cref{punc-thm}. 

\bthm[\S\ref{main-thm-pf}] \lab{main-thm}
For a Henselian, regular, local ring $R$ of dimension $\ge 2$ whose residue field is of dimension $\le 1$ \up{in the sense recalled in Definition \uref{dim-1-def}} and an $R$-torus $T$, we have
\[
H^1(U_R, T) = H^2(U_R, T) = 0.
\]
\ethm

The vanishing of $H^1(U_R, T)$, included here for completeness, follows already from \cite{CTS79}*{6.9}.

\bpp[Notation and conventions] \lab{conv}
For a semilocal ring $R$, we let 
\[
U_R \subset \Spec(R)
\]
 be the open complement of the closed points. For most schemes $S$ that we consider, we have 
 \[
 H^2(S, \bG_m) = H^2(S, \bG_m)_\tors
 \]
 (see \Cref{reg-inj}), so we phrase our results about the (cohomological) Brauer group in terms of \'{e}tale cohomology. Other than in the proof of \Cref{GR-inj}, we do not use the relationship with Azumaya algebras. For a scheme morphism $S' \ra S$, we let $(-)_{S'}$ denote base change to $S'$ and let $\Res_{S'/S}(-)$ denote the restriction of scalars. For a field $k$, we let $\ov{k}$ denote a fixed choice of its algebraic closure. We let $W(-)$ denote the $p$-typical Witt vectors. When no confusion seems likely, we let $\cO$ abbreviate the structure sheaf $\cO_S$ of a scheme $S$. 
\epp

}

\subsection*{Acknowledgements}
\ready{I thank Peter Scholze for several very helpful conversations. I thank Jean-Louis Colliot-Th\'{e}l\`{e}ne, Ofer Gabber, Ariyan Javanpeykar, Daniel Loughran, Bjorn Poonen, and Minseon Shin for helpful conversations and correspondence. I thank the referees for very thorough and helpful comments and suggestions.
}



\section{Passage to a finite flat cover} \lab{finite-flat}

\ready{
The perfectoid approach to the purity conjecture hinges on the ability to pass to an infinitely ramified cover of a regular local ring $R$ without killing Brauer classes of its punctured spectrum. The results of the present section facilitate this. To highlight the inputs to their proofs, we chose an axiomatic approach when presenting the key \Cref{pump-R-0,pump-R}. Concrete situations in which these propositions apply are described in \Cref{cor-H2-et,cor-H2-inj} and Remark \ref{rem-lci}.

\blem \lab{rep-prop}
For a finite, locally free scheme morphism $\pi\colon S' \ra S$ and an $S$-affine $S$-group scheme $G$, the homogeneous space 
\[
X \ce (\Res_{S'/S}(G_{S'}))/G
\]
is representable by an $S$-affine scheme that is smooth if so is $G$. In addition, if $\pi$ has a section, then
\be \lab{split}
\Res_{S'/S}(G_{S'}) \cong G \times_S X
\ee
as $S$-schemes and, in the case when $G$ is commutative, even as $S$-group schemes.
\elem

\bpf
Both the representability by an $S$-affine scheme and the smoothness are properties that are fppf local on $S$, so, by base change along $\pi$, we assume that $\pi$ has a section:
\[
\xymatrix{S \ar[r]^-{s} \ar@{=}[rd] & S' \ar[d]^-{\pi}  \\ & S. }
\]
Then the adjunction map 
\[
i\colon G \hra \Res_{S'/S}(G_{S'}) \qq \x{has a section} \qq j \colon \Res_{S'/S}(G_{S'}) \surjects \Res_{S/S}(G) \cong G,
\]
which is a group morphism. It follows that $X \cong \Ker j$ over $S$, compatibly with group structures if $G$ is commutative. Since $\Res_{S'/S}(G_{S'})$ is an $S$-affine $S$-group scheme (see \cite{BLR90}*{7.6/4 and its proof}), the representability of $X$ by an $S$-affine scheme and the decomposition \eqref{split} follow. If $G$ is smooth, then so is $\Res_{S'/S}(G_{S'})$, and hence $X$ is, too (see \cite{BLR90}*{7.6/5} and \cite{SGA3Inew}*{VI$_{\text{B}}$,~9.2~(xii)}).
\epf

\bprop \lab{pump-R-0}
For a finite, flat map $R \ra R'$ of local rings, an open subscheme $V \subset \Spec R$, and an affine, smooth $R$-group scheme $G$, if
\benuma
\item \lab{PR0-1}
$\Gamma(\Spec R, \cO) \cong \Gamma(V, \cO)$ via pullback\uscolon and 

\item \lab{PR0-2}
every $G$-torsor is trivial over $R$\uscolon 
\eenum
then the following pullback is injective\ucolon
\be \lab{PR0-eq}
H^1_\et(V, G) \hra H^1_\et(V_{R'}, G).
\ee
\eprop

\bpf
By \Cref{rep-prop}, the homogeneous space 
\[
X \ce (\Res_{R'/R}(G_{R'}))/G
\]
is an affine $R$-scheme. Thus, due to \ref{PR0-1}, we have 
\[
X(R) \cong X(V) \q \x{via pullback.}
\]
However, due to \ref{PR0-2}, every element of $X(R)$ lifts to $(\Res_{R'/R}(G_{R'}))(R)$. Consequently, every element of $X(V)$ lifts to $(\Res_{R'/R}(G_{R'}))(V)$, so, by \cite{Gir71}*{III.3.2.2}, the map
\be \lab{PR0-eq2}
H^1_\et(V, G) \ra H^1_\et(V, \Res_{R'/R}(G_{R'}))
\ee
is injective. However, the projection $\pi\colon V_{R'} \ra V$ is finite, so, as may be checked on strict Henselizations at points of $V$, the \'{e}tale sheaf $R^1 \pi_*(G_{V_{R'}})$ vanishes. By \cite{Gir71}*{V.3.1.3}, this implies that 
\[
H^1_\et(V, \Res_{R'/R}(G_{R'})) \cong H^1_\et(V_{R'}, G),
\]
so the injectivity of \eqref{PR0-eq} follows from that of \eqref{PR0-eq2}.
\epf

\bprop \lab{pump-R}
For a finite, flat map $R \ra R'$ of local rings, an open subscheme $V \subset \Spec R$, and an $R$-torus $T$ that splits over $R'$, if
\benuma
\item \lab{PR-1}
$\Gamma(\Spec R, \cO) \cong \Gamma(V, \cO)$ via pullback\uscolon

\item \lab{PR-2}
every $(\Res_{R'/R}(T_{R'}))/T$-torsor is trivial over $R$\uscolon and

\item \lab{PR-3}
$\Pic(V_{R' \tensor_R R'}) = 0$\uscolon
\eenum
then the following pullback is injective\ucolon
\be \lab{pump-R-eq-0}
H^2_\et(V, T) \hra  H^2_\et(V_{R'} , T);
\ee
if instead of \ref{PR-3} we have
\begin{enumerate}[label={{\upshape(\arabic*${}^\prime$)}}] \addtocounter{enumi}{2}
\item \lab{PR-3-pr}
$\Pic(V_{R' \tensor_R R'})$ is torsion free and $\Pic(V_{R'})$ is torsion\uscolon
\eenum
\up{in addition to \ref{PR-1} and \ref{PR-2}}, then the pullback is injective on the torsion subgroups\ucolon
\be \lab{pump-R-eq}
H^2_\et(V, T)_\tors \hra  H^2_\et(V_{R'} , T)_\tors.
\ee
\eprop

\bpf
By \Cref{rep-prop}, the quotient $G \ce (\Res_{R'/R}(T_{R'}))/T$ is representable by an affine, smooth $R$-group scheme and $G_{R'}$ is a direct factor $R'$-group scheme of $(\Res_{(R' \tensor_R R')/R'}(\bG_m))^{\rk T}$. In particular, \ref{PR-1}--\ref{PR-2} ensure that \Cref{pump-R-0} applies to $G$, and we conclude the injectivity of the maps
\[
\tst H^1_\et(V, G) \hra H^1_\et(V_{R'}, G) \hra \bigoplus_{i = 1}^{\rk T} H^1_\et(V_{R'}, \Res_{(R' \tensor_R R')/R'}(\bG_m)) \cong (\Pic(V_{R' \tensor_R R'}))^{\rk T},
\]
where the identification follows from the exactness in the \'{e}tale topology of the pushforward along a finite morphism  (see \cite{SGA4II}*{VIII, 5.5}). The assumption \ref{PR-3} then gives $H^1_\et(V, G) = 0$ and hence, due to the cohomology sequence
\be \lab{coho-LES}
\dotsc \ra H^1_\et(V_{R'}, T) \ra H^1_\et(V, G) \ra H^2_\et(V, T) \ra H^2_\et(V_{R'}, T) \ra \dotsc,
\ee
implies the claimed \eqref{pump-R-eq-0}. In addition, since $T$ splits over $R'$, we have 
\[
H^1_\et(V_{R'}, T) \cong (\Pic(V_{R'}))^{\rk T}.
\]
Thus, if \ref{PR-3-pr} holds instead, then $H^1_\et(V, G)$ is torsion free and injects into $H^2_\et(V, T)$, to the effect that then \eqref{coho-LES} implies \eqref{pump-R-eq}.
\epf

\bcor \lab{cor-H2-et}
For a Henselian, regular, local ring $R$ of dimension $\ge 2$ whose residue field $k$ is of dimension $\le 1$, a nonzero finite \'{e}tale $R$-algebra $R'$, and an $R$-torus $T$, the following pullbacks are injective\ucolon
\[
H^1_\et(U_R, T) \hra H^1_\et(U_{R'}, T) \qq \text{and} \qq H^2_\et(U_R, T) \hra H^2_\et(U_{R'}, T).
\]
\ecor

\bpf
We set $V \ce U_R$, so that, due to the $R$-finiteness of $R'$, we have $V_{R'} = U_{R'}$. We lose no generality by enlarging $R'$, so we assume that $R'$ is local and $T$ splits over $R'$. Thus, the claim follows from \Cref{pump-R-0,pump-R} once we explain why their assumptions \ref{PR0-1}--\ref{PR0-2} and \ref{PR-1}--\ref{PR-3}~hold.

The assumption \ref{PR-1} holds because $R$ is Noetherian of depth $\ge 2$ (see \cite{EGAIV2}*{5.10.5}). Since $R'/R$ is finite \'{e}tale, $(\Res_{R'/R}(T_{R'}))/T$ and $T$ are both tori, so, by \Cref{dim-1-tori}, their torsors are trivial over $k$. Thus, since $R$ is Henselian, the same is true over $R$ (see \cite{EGAIV4}*{18.5.17}), so \ref{PR-2} holds. Finally, \ref{PR-3} holds because $R' \tensor_R R'$ is a product of regular local rings of dimension $\ge 2$.
\epf

\bcor \lab{cor-H2-inj}
For a finite, flat map $f\colon R \ra R'$ of strictly Henselian, regular, local rings of dimension $\ge 4$, the following pullback is injective\ucolon
\[
H^2_\et(U_R, \bG_m) \hra H^2_\et(U_{R'}, \bG_m).
\]
\ecor

\bpf
We apply \Cref{pump-R} with $V = U_R$. Its assumption \ref{PR-1} holds because $R$ is of depth $\ge 2$. Since $R$ is strictly Henselian and $(\Res_{R'/R}((\bG_m)_{R'}))/\bG_m$ is $R$-smooth and affine (see \Cref{rep-prop}), the assumption \ref{PR-2} holds, too. Since $R$ and $R'$ are regular, $f$ is necessarily a local complete intersection morphism (see \cite{SP}*{\href{http://stacks.math.columbia.edu/tag/0E9K}{0E9K}}), so the local ring $R' \tensor_R R'$ is a local complete intersection of dimension $\ge 4$ (see \cite{SP}*{\href{http://stacks.math.columbia.edu/tag/069I}{069I}, \href{http://stacks.math.columbia.edu/tag/07D3}{07D3}, \href{http://stacks.math.columbia.edu/tag/09Q7}{09Q7}}). However, by the Grothendieck--Lefschetz theorem \cite{SGA2new}*{XI, 3.13 (ii)}, the Picard group of the punctured spectrum of a local complete intersection of dimension $\ge 4$ vanishes, so \ref{PR-3} holds.
\epf

\brems
\remi \lab{rem-lci}
As is clear from its proof, one may strengthen \Cref{cor-H2-inj} by assuming instead that $f$ is a finite, flat, local complete intersection morphism of strictly Henselian, Noetherian, local rings that are local complete intersections of dimension $\ge 4$. 


\remi
Gabber pointed out the following more direct argument for \Cref{cor-H2-inj}. Consider the spectral sequence \cite{SGA4II}*{V.3.3} for the finite flat covering $U_{R'}/U_R$ and $\bG_m$ as coefficients:
\be \lab{Gab-ss}
\qq E_2^{ij} = H^i(U_{R'}/U_R, H^j_\et(-, \bG_m)) \Ra H^{i + j}_\et(U_R, \bG_m),
\ee
where we use \cite{Gro68c}*{11.7} to identify the fppf and the \'{e}tale cohomologies with $\bG_m$ coefficients. Each self-product $R' \tensor_R \dotsc \tensor_R R'$ is a strictly Henselian local complete intersection of dimension $\ge 4$ (compare with the proof of \Cref{cor-H2-inj}) and its punctured spectrum is $U_{R'} \times_{U_R} \dotsc \times_{U_R} U_{R'}$, so, by the Grothendieck--Lefschetz theorem,  $E_2^{i1} = 0$ for every $i$. Moreover, each $R' \tensor_R \dotsc \tensor_R R'$ has depth $\ge 2$, so the $E_2^{i0}$ terms are identified with the corresponding terms of the analogous spectral sequence for the covering $R'/R$, to the effect that $E_2^{i0} = 0$ for $i > 0$. Thus, the $i + j = 2$ diagonal of the $E_2$ page of \eqref{Gab-ss} gives the desired injectivity
\[
\qq H^2(U_R, \bG_m) \hra E_2^{02} \subset H^2(U_{R'}, \bG_m).
\]
\erems
}



\section{Passage to the completion} \lab{completion}

\ready{
We will need the flexibility of replacing a Henselian ring by its completion without killing Brauer classes.
The following standard results achieve this. Their general theme goes back at least to \cite{Elk73} and they rely on the work of Gabber \cite{Gab81} and Gabber--Ramero \cite{GR03}.

\blem \lab{GR-inj}
For a ring $R$ that is Henselian along a principal ideal $(f) \subset R$ generated by a nonzerodivisor $f\in R$,  the following pullback, where $\wh{R}$ denotes the $f$-adic completion of $R$, is bijective\ucolon
\be \lab{GR-inj-eq}
\tst H^2_\et(R[\f{1}{f}], \bG_m)_\tors \isomto H^2_\et(\wh{R}[\f{1}{f}], \bG_m)_\tors.
\ee
\elem

\bpf
By \cite{GR03}*{5.4.41}, for every $n \ge 0$, the following pullback is bijective:
\[
\tst H^1_\et(R[\f{1}{f}], \GL_n) \isomto H^1_\et(\wh{R}[\f{1}{f}], \GL_n).
\]
In addition, for any two $(\PGL_n)_{R[\f{1}{f}]}$-torsors $X$ and $X'$, their isomorphism functor $\Iso_{\PGL_n}(X, X')$ is representable by an affine, smooth  $R[\f{1}{f}]$-scheme (which \'{e}tale locally on $R[\f{1}{f}]$ is isomorphic to $(\PGL_n)_{R[\f{1}{f}]}$). Consequently, by \cite{GR03}*{5.4.21}, if $X$ and $X'$ are not isomorphic, they cannot become isomorphic over $\wh{R}[\f{1}{f}]$, to the effect that the following pullback map is injective:
\[
\tst H^1_\et(R[\f{1}{f}], \PGL_n) \hra H^1_\et(\wh{R}[\f{1}{f}], \PGL_n).
\]
Thus, the nonabelian cohomology exact sequences of \cite{Gir71}*{IV.4.2.10} that result from the central extension $1 \ra \bG_m \ra \GL_n \ra \PGL_n \ra 1$ fit into the commutative diagram
\be\ba \lab{GL-PGL-diag}
\xymatrix{
\tst H^1_\et(R[\f{1}{f}], \GL_n) \ar[d]^{\wr} \ar[r] & H^1_\et(R[\f{1}{f}], \PGL_n) \ar[r] \ar@{^(->}[d] & H^2_\et(R[\f{1}{f}], \bG_m) \ar[d] \\
    H^1_\et(\wh{R}[\f{1}{f}], \GL_n) \ar[r] & H^1_\et(\wh{R}[\f{1}{f}], \PGL_n) \ar[r] & H^2_\et(\wh{R}[\f{1}{f}], \bG_m).
}
\ea \ee
This diagram shows that no nonzero element of the image of $H^1_\et(R[\f{1}{f}], \PGL_n)$ in $H^2_\et(R[\f{1}{f}], \bG_m)$ maps to zero in $H^2_\et(\wh{R}[\f{1}{f}], \bG_m)$. By \cite{Gab81}*{II, Thm.~1}, as $n$ varies, these images sweep out $H^2_\et(R[\f{1}{f}], \bG_m)_\tors$, so the injectivity \eqref{GR-inj-eq} follows. 

By also applying \emph{loc.~cit.} to $\wh{R}[\f{1}{f}]$ and using \cite{GR03}*{5.8.14},\footnote{The assumptions of \cite{GR03}*{5.8.14} are met for $G = \PGL_n$ by \cite{GR03}*{5.8.5}; in fact, they are met whenever $G$ is semisimple because, by \cite{Tho87}*{3.2~(3)}, over an affine base such a $G$ embeds into some $\GL_N$ and the quotient $\GL_N/G$ is affine (see \cite{Alp14}*{9.4.1 and 9.7.5}) so its structure sheaf $\sO_{\GL_N/G}$ is ample and $\GL_N$-equivariant.} which ensures that the middle vertical arrow in \eqref{GL-PGL-diag} is bijective, we conclude that the map \eqref{GR-inj-eq} is also surjective. (We will only use the injectivity of \eqref{GR-inj-eq}, so we focused the proof on this simpler part of the claim.) 
\epf

To deduce \Cref{H2-inj} from \Cref{GR-inj}, we will use the following widely-known result.

\blem[\cite{Gro68b}*{1.8}, see also Remark \ref{flasque}] \lab{reg-inj}
For a Noetherian, integral, regular scheme $X$ and its function field $K$, the pullback
\[
H^2_\et(X, \bG_m) \hra H^2_\et(K, \bG_m)
\]
is injective\uscolon in particular, $H^2_\et(X, \bG_m)$ is torsion. \QED
\elem

\bprop \lab{H2-inj}
For a Henselian, regular, local ring $(R, \fm)$, the following pullback, where $\wh{R}$ denotes the $\fm$-adic completion of $R$, is injective\ucolon
\be \lab{H2-inj-eq}
H^2_\et(U_R, \bG_m) \hra H^2(U_{\wh{R}}, \bG_m).
\ee
\eprop

\bpf
Let $f_1, \dotsc, f_{\dim(R)} \in \fm$ be a regular sequence that generates $\fm$, set $R_0 \ce R$, and, for each $1 \le i \le \dim(R)$, let $R_i$ be the $f_i$-adic completion of $R_{i - 1}$. Explicitly, each $R_i$ is the $(f_1, \dotsc, f_i)$-adic completion of $R$: indeed, by induction on $i$, this follows by forming $\varprojlim_n$ of the short exact sequences
\[
0 \ra R/(f^n_1, \dotsc, f_{i - 1}^n) \xra{f^m_i} R/(f^n_1, \dotsc, f_{i - 1}^n) \ra R/(f^n_1, \dotsc, f_{i - 1}^n, f_i^m) \ra 0 \q \text{for} \q i > 1, \q m \ge 1.
\]
In particular, $R_{\dim(R)} \cong \wh{R}$ and each $R_i$ is local, regular, and Henselian (see \cite{SP}*{\href{http://stacks.math.columbia.edu/tag/0AGX}{0AGX}, \href{http://stacks.math.columbia.edu/tag/07NY}{07NY}, \href{http://stacks.math.columbia.edu/tag/0DYD}{0DYD}}). Consequently, for $1 \le i \le \dim(R)$, \Cref{GR-inj,reg-inj} give the commutative diagram
\[
\xymatrix{
H^2_\et(U_{R_{i - 1}}, \bG_m) \ar@{^(->}[d]\ar[r] & H^2_\et(U_{R_{i}}, \bG_m) \ar@{^(->}[d] \\
H^2_\et(R_{i - 1}[\f{1}{f_{i}}], \bG_m) \ar@{^(->}[r] & H^2_\et(R_{i}[\f{1}{f_{i}}], \bG_m),
}
\]
which shows the injectivity of its top horizontal map. Induction on $i$ then gives \eqref{H2-inj-eq}.
\epf

}



\section{The $p$-primary Brauer group in the perfectoid case} \lab{perf-section}

\ready{
While \S\S\ref{finite-flat}--\ref{completion} facilitate passage to perfect or perfectoid rings, the present one investigates the $p$-primary part of the Brauer group of such a ring. We begin with the simpler positive characteristic~case.

\bprop \lab{perf-Gm}
For a prime $p$ and a perfect $\bF_p$-scheme $X$, we have
\[
H^i_\et(X, \bG_m)[p^\infty] = 0 \q \text{for every} \q i \in \bZ.
\]
\eprop

\bpf
Every \'{e}tale $X$-scheme inherits perfectness from $X$ (see \cite{SGA5}*{XV, \S1, Prop.~2 c) 2)}). Therefore, on the \'{e}tale site of $X$, the $p$-power map is an automorphism of the sheaf $\bG_m$.
\epf

A mixed characteristic analogue of \Cref{perf-Gm} is \Cref{p-oid-Br} below, which concerns perfectoid rings.  The latter were introduced by Scholze in \cite{Sch12} in the context of rigid geometry, with variants in other contexts appearing afterwards. Their axiomatics that suit our purposes are captured by the following definition and discussion, which  are related to \cite{BMS16}*{\S3.2}.

\bd \lab{perf-def}
For a prime $p$, a $p$-torsion free ring $R$ is \emph{perfectoid} if $R$ is $p$-adically complete and the divisor $(p) \subset R$ has a $p$-th root in the sense that there is a $\varpi \in R$ with $(\varpi^p) = (p)$, and 
\be \lab{varpi-def}
\xymatrix{R/\varpi \ar[r]_-{\sim}^{x\, \mapsto\, x^p} & R/p.}
\ee
(Since $(\varpi) \subset R$ is the preimage of the kernel of the Frobenius of $R/p$, it is uniquely~determined.)
\ed

\brems
\remi
The $p$-torsion freeness, the $p$-adic completeness, and \eqref{varpi-def} imply that $R$ is reduced.

\remi \lab{pn-roots}
The $p$-adic completeness of $R$ implies that the reduction modulo $p$ map
\[
\tst \qq \varprojlim_{x\mapsto x^p} R \isomto \varprojlim_{x\mapsto x^p} (R/p)
\]
is an isomorphism of multiplicative monoids (see the proof of \cite{Sch12}*{3.4~(i)}). Thus, due to the surjectivity of \eqref{varpi-def}, there is a $p$-power compatible sequence $(\dotsc, \varpi_2, \varpi_1)$ of elements of $R$ with $\varpi_1 \equiv \varpi \bmod p$. Since $(p) = (\varpi^p)$, this gives 
\[
\qqq (\varpi) = (\varpi_1) \subset (\varpi_2) \subset \dotsc  \qq \text{and} \qq (\varpi_n^{p^n}) = (p) \q \text{for every} \q n > 0.
\]
In particular, each $(\varpi_n) \subset R$ is uniquely determined: by induction on $n$, it is the preimage of the kernel of the $p^n$-power Frobenius of $R/p$, so that
\be \lab{pn-root}
\qq \xymatrix{R/\varpi_n \ar[r]_-{\sim}^-{x \,\mapsto\, x^{p^n}} & R/p.}
\ee

\remi \lab{more-surj}
By \eqref{varpi-def}, modulo $p^2$ every element of $R$ is of the form $x^p + p y^p$ or, equivalently, $x^p + \varpi^p y^{\prime p}$. In particular, modulo $p\varpi$ every element of $R$ is a $p$-th power (a special case of \cite{BMS16}*{3.9}).

\remi
By \cite{BMS16}*{3.10~(ii)}, an $R$ as in \Cref{perf-def} is perfectoid in the sense of the definition \cite{BMS16}*{3.5}. Conversely, a $p$-torsion free ring that is perfectoid in the sense of \emph{loc.~cit.}~is perfectoid in the sense of \Cref{perf-def} due to \cite{BMS16}*{3.9 and 3.10~(i)}. Compatibility with other possible definitions is discussed in \cite{BMS16}*{\S3.2, especially 3.20}.
\erems

The following simple lemma often helps to recognize perfectoid rings in nature.

\blem \lab{p-oid-rec}
For a prime $p$ and a $p$-torsion free ring $R$ such that $(\varpi^p) = (p)$ for some $\varpi \in R$, if $R$ is integrally closed in $R[\f{1}{p}]$, then the map $R/\varpi \xra{x\, \mapsto\, x^p} R/p$ is injective\uscolon if, in addition, every element of $R/p$ is a $p$-th power, then the $p$-adic completion $\wh{R}$ of $R$ is perfectoid.
\elem

\bpf
If $r \in R$ represents a class in the kernel of the map 
\[
\tst R/\varpi \xra{x\, \mapsto\, x^p} R/p,
\]
then $r^p = \varpi^ps$ for some $s \in R$. Since $R[\f{1}{\varpi}] = R[\f{1}{p}]$ and $R$ is integrally closed in $R[\f{1}{p}]$, this implies that $\f{r}{\varpi} \in R$. Thus, $r \in (\varpi)$, and the injectivity follows. Since $\wh{R}/\varpi \cong R/\varpi$ and $\wh{R}/p \cong R/p$, the second assertion follows as well.
\epf

To study the $p$-primary Brauer group of $R[\f{1}{p}]$, we will use the tilting equivalence of Scholze \cite{Sch12}*{7.12}. More precisely, since we do not wish to restrict to $R[\f{1}{p}]$ that are algebras over some perfectoid field, we will use the version of this equivalence presented by Kedlaya and Liu in \cite{KL15}. We will review its precise statement in \S\ref{tilt}, after discussing the following auxiliary reduction.

\bpp[A reduction to the case $R = (R\lb\f{1}{p}\rb)^\circ$] \lab{powerbdd}
We endow a $p$-torsion free perfectoid ring $R$ with its $p$-adic topology and $R[\f{1}{p}]$ with the unique ring topology for which $R \subset R[\f{1}{p}]$ is open, so that $R[\f{1}{p}]$ is a Tate ring in the sense of Huber (see \cite{Hub93a}). Due to \eqref{pn-root}, if $x \in R$ is such that $x^{p^n} \in p^{p^n}R$, then $x \in pR$; in particular, $R$ contains the topologically nilpotent elements $(R[\f{1}{p}])^{\circ \circ}$. Thus, since the subring $(R[\f{1}{p}])^\circ \subset R[\f{1}{p}]$ of powerbounded elements is the union of the open, bounded subrings of $R[\f{1}{p}]$ that contain $R$ (see \cite{Hub93a}*{1.2--1.3}), we conclude that, in the notation of Remark~\ref{pn-roots},
\be \lab{almost-cok}
\tst \text{each} \q \varpi_n \q \text{kills the cokernel of the inclusion} \q R \subset (R[\f{1}{p}])^\circ
\ee
(a special case of \cite{BMS16}*{3.21}). In particular, the subring $(R[\f{1}{p}])^\circ \subset R[\f{1}{p}]$ is bounded (that is, $R[\f{1}{p}]$ is \emph{uniform}), so $(R[\f{1}{p}])^\circ$ is $p$-adically complete. In fact, $(R[\f{1}{p}])^\circ$ is even perfectoid: indeed, the~map
\[
\tst (R[\f{1}{p}])^\circ/\varpi \xra{x\,\mapsto\, x^p} (R[\f{1}{p}])^\circ/p
\]
is injective because $x^p = \varpi^p y$ in $(R[\f{1}{p}])^\circ$ implies $\f{x}{\varpi} \in (R[\f{1}{p}])^\circ$; 
it is also surjective because, by \eqref{almost-cok} and Remark \ref{more-surj}, for every $x \in (R[\f{1}{p}])^\circ$ we have $\varpi_1 x = y^p + p\varpi_1 z$ with $y, z \in R$, so that $\f{y}{\varpi_2} \in (R[\f{1}{p}])^\circ$.

In conclusion, by replacing $R$ by $(R[\f{1}{p}])^\circ$, we reduce the study of $R[\f{1}{p}]$ to the case when $R = (R[\f{1}{p}])^\circ$. Then $R$ is a ring of integral elements, so that $(R[\f{1}{p}], R)$ is an affinoid Tate ring (see \cite{Hub93a}*{\S3}).
\epp

\bpp[The tilting equivalence] \lab{tilt}
Let $R$ be a $p$-torsion free perfectoid ring with $R = (R[\f{1}{p}])^\circ$. The norm function 
\be \lab{norm}
x\mapsto \inf (\{ 2^{n} \, \vert\, n \in \bZ \text{ with } p^{n} x \in R\})
\ee
makes $R[\f{1}{p}]$ a Banach $\bQ_p$-algebra (in the sense of \cite{KL15}*{2.2.1}) whose unit ball is $R$. In particular, the pair $(R[\f{1}{p}], R)$ becomes a \emph{perfectoid Banach $\bQ_p$-algebra} in the sense of \cite{KL15}*{3.6.1} (see \cite{KL15}*{3.6.2~(e)}). Its \emph{tilt} is
\[
\tst (R^\flat[\f{1}{\varpi^\flat}], R^\flat), \qq \text{where} \qq R^\flat \ce \varprojlim_{x \mapsto x^p} (R/p) \q \text{and} \q \varpi^\flat \overset{\ref{pn-roots}}{\ce} (\dotsc, \varpi_2 \bmod p, \varpi_1 \bmod p) \in R^\flat,
\]
so that $\varpi^\flat \in R^\flat$ is a nonzerodivisor and $R^\flat$ is a $\varpi^\flat$-adically complete, perfect $\bF_p$-algebra. We endow $R^\flat$ with its $\varpi^\flat$-adic topology and $R^\flat[\f{1}{\varpi^\flat}]$ with the unique ring topology for which $R^\flat \subset R^\flat[\f{1}{\varpi^\flat}]$ is open. Due to the compatible multiplicative monoid isomorphisms
\[
\tst R^\flat \overset{\ref{pn-roots}}{\cong} \varprojlim_{x \mapsto x^p} R \qq \text{and} \qq R^\flat[\f{1}{\varpi^\flat}] \cong \varprojlim_{x \mapsto x^p} (R[\f{1}{p}]),
\]
$R^\flat = (R^\flat[\f{1}{\varpi^\flat}])^\circ$. The norm \eqref{norm} with $\varpi^\flat$ in place of $p$ makes $(R^\flat[\f{1}{\varpi^\flat}], R^\flat)$ a Banach $\bF_p$-algebra. 

By \cite{KL15}*{3.1.13, 3.6.15}, the structure presheaves of the spaces
\be \lab{Spa-Spa}
\tst \Spa(R[\f{1}{p}], R) \qq \text{and} \qq \Spa(R^\flat[\f{1}{\varpi^\flat}], R^\flat)
\ee
are sheaves, that is, these spaces are adic. Moreover, by \cite{KL15}*{3.6.14}, 
the two spaces in \eqref{Spa-Spa} are naturally (and functorially in $R$) homeomorphic in such a way that rational subsets correspond to rational subsets. In addition, by the almost purity theorem in this context \cite{KL15}*{3.6.23}, this homeomorphism extends to an equivalence of \'{e}tale sites\footnote{\lab{foot}\ready{The \'{e}tale sites are defined as in \cite{Sch12}*{7.1, 7.11}: e.g.,~a morphism to $\Spa(R[\f{1}{p}], R)$ is \'{e}tale if and only if on an open cover of the source it is an open immersion followed by a finite \'{e}tale map to a rational subspace of $\Spa(R[\f{1}{p}], R)$; for stability of \'{e}tale morphisms under compositions and fiber products, see \cite{KL15}*{8.2.17 (c)}.}} 
\be \lab{et-sit-eq}
\tst \Spa(R[\f{1}{p}], R)_\et \cong \Spa(R^\flat[\f{1}{\varpi^\flat}], R^\flat)_\et
\ee
that identifies finite \'{e}tale $(R[\f{1}{p}])$-algebras and finite \'{e}tale $(R^\flat[\f{1}{\varpi^\flat}])$-algebras. 
\epp


\bthm \lab{p-oid-Br}
For a $p$-torsion free perfectoid ring $R$ and a commutative, finite, \'{e}tale $R[\f{1}{p}]$-group scheme $G$ of $p$-power order, we have 
\be\lab{POB-eq}
\tst H^i_\et(R[\f{1}{p}], G) = 0 \q \text{for} \q i \ge 2, \qq \text{so also} \qq H^i_\et(R[\f{1}{p}], \bG_m)[p^\infty] = 0 \q \text{for} \q i \ge 2.
\ee
\ethm

\bpf
\ready{
By \S\ref{powerbdd}, we may assume that $R = (R[\f{1}{p}])^\circ$. Then \cite{Hub96}*{3.2.9} (granted that we explain why it applies, as we do below; we choose $U = \Spa A$ in \emph{loc.~cit.}) 
gives the identification
\be \lab{iso-1}
\tst H^i_\et(R[\f{1}{p}], G) \cong H^i_\et(\Spa(R[\f{1}{p}], R), G).
\ee
Since the equivalence \eqref{et-sit-eq} identifies finite \'{e}tale $(R[\f{1}{p}])$-algebras and finite \'{e}tale $(R^\flat[\f{1}{\varpi^\flat}])$-algebras, $G$ determines a commutative, finite, \'{e}tale $(R^\flat[\f{1}{\varpi^\flat}])$-group scheme $G^\flat$ of $p$-power order such that
\be \lab{iso-2}
\tst H^i_\et(\Spa(R[\f{1}{p}], R), G) \cong H^i_\et(\Spa(R^\flat[\f{1}{\varpi^\flat}], R^\flat), G^\flat).
\ee
By \cite{Hub96}*{3.2.9} again (with the same caveat), 
\be \lab{iso-3}
\tst H^i_\et(\Spa(R^\flat[\f{1}{\varpi^\flat}], R^\flat), G^\flat) \cong H^i_\et(R^\flat[\f{1}{\varpi^\flat}], G^\flat).
\ee
However, by \cite{SGA4III}*{X, 5.1}, the \'{e}tale cohomological $p$-dimension of an affine Noetherian $\bF_p$-scheme is at most $1$, so, by a limit argument, 
\[
\tst H^i_\et(R^\flat[\f{1}{\varpi^\flat}], G^\flat) = 0 \q \x{for} \q i \ge 2.
\]
Due to \eqref{iso-1}--\eqref{iso-3}, this gives the desired 
\[
\tst H^i_\et(R[\f{1}{p}], G) = 0 \q \x{for}\q i \ge 2.
\]
The second part of \eqref{POB-eq} follows by choosing $G = \mu_p$.

In order to ensure that the structure presheaves of adic spectra are sheaves, the book \cite{Hub96} is written under a blanket Noetherianness assumption \cite{Hub96}*{1.1.1}. Thus, the deduction of \eqref{iso-1} and \eqref{iso-3} above from \cite{Hub96}*{3.2.9} implicitly involves the following limit argument. 

The ring $R$ is a filtered direct limit of $p$-torsion free $\bZ_p$-subalgebras $R_j$ of finite type that are integrally closed in $R_j[\f{1}{p}]$ (to ensure the latter, we use the reducedness of $R$ and \cite{EGAIV2}*{7.8.6~(ii), 7.8.3~(ii)--(iii)}; the purpose of this integral closedness is to be able to later write $\Spa(R_j[\f{1}{p}], R_j)$). In fact, since $R$ is $p$-adically complete, it is even the filtered direct limit of the Henselizations $R_j^h$ of $R_j$ along the ideals $(p) \subset R_j$. We may assume that $G$ descends to each $R_j[\f{1}{p}]$, so, by \cite{SGA4II}*{VII,~5.9},
\be \lab{Hi-alg-lim}
\tst H^i_\et(R[\f{1}{p}], G) \cong \varinjlim_j H^i_\et(R_j^h[\f{1}{p}], G) \qq \text{for every $i$.}
\ee
}
We equip each $R_j$ with the $p$-adic topology and $R_j[\f{1}{p}]$ with the unique ring topology for which $R_j \subset R_j[\f{1}{p}]$ is open. Then a valuation of $R_j[\f{1}{p}]$ whose values on $R_j$ are $\le 1$ is continuous if and only if the values of $\{p^n\}_{n > 0}$ are not bounded below, and likewise for $R[\f{1}{p}]$. In particular, the map
\be \lab{spa-homeo}
\tst \Spa(R[\f{1}{p}], R) \ra \varprojlim_j (\Spa(R_j[\f{1}{p}], R_j))
\ee
is a homeomorphism that respects rational subsets. In fact, by following the arguments of \cite{Sch17}*{proof of 6.4~(ii)},\footnote{More precisely, the argument of \emph{loc.~cit.}~shows the following. Let $\{(A_j[\f{1}{\varpi}], A_j)\}_{j \in J}$ be a filtered direct system of affinoid Tate rings that have a common pseudouniformizer $\varpi$ and such that each $A_j$ is Noetherian, a ring of definition for $A_j[\f{1}{\varpi}]$, and integrally closed in $A_j[\f{1}{\varpi}]$ (so that the f-adic rings $A_j[\f{1}{\varpi}]$ meet the assumption \cite{Hub96}*{1.1.1 a)} and the results of \emph{op.~cit.}~apply). Suppose that the $\varpi$-adically completed direct limit $A \ce \p{\varinjlim_{j \in J} A_j}\wh{\ }$ has $A[\f{1}{\varpi}]$ to be perfectoid (in the sense of \cite{Sch17}*{3.1}; see \cite{BMS16}*{3.20} for the compatibility with \Cref{perf-def}). Then the base change functors induce an equivalence
\[
\tst 2\x{-}\varinjlim_{j \in J}\p{ \Spa(A_j[\f{1}{\varpi}], A_j)_{\et,\, \mathrm{qcqs}}} \isomto \Spa(A[\f{1}{\varpi}], A)_{\et,\, \mathrm{qcqs}}
\]
that relates the full subcategories consisting of quasi-compact and quasi-separated objects of the indicated \'{e}tale sites (defined as in \cref{foot}). The (irrelevant for the proof) difference from \cite{Sch17}*{6.4~(ii)} is that there one assumes each $(A_j[\f{1}{\varpi}], A_j)$ to be affinoid perfectoid in the sense of \emph{op.~cit.}} we see that \eqref{spa-homeo} extends to an equivalence between the \'{e}tale site $\Spa(R[\f{1}{p}], R)_\et$ and the $2$-limit of the \'{e}tale sites $\Spa(R_j[\f{1}{p}], R_j)_\et$ granted that we restrict to quasi-compact and quasi-separated adic spaces in these sites (which does not change the associated topoi). 
As in the proof of \cite{SGA4II}*{VII, 5.7}, generalities on projective limits of fibered topoi then imply\footnote{\ready{Alternatively and more concretely, one may use hypercoverings and \cite{SP}*{\href{http://stacks.math.columbia.edu/tag/01H0}{01H0}} to deduce \eqref{Hi-an-lim}.}} that 
\be \lab{Hi-an-lim}
\tst H^i_\et(\Spa(R[\f{1}{p}], R), G) \cong \varinjlim_j H^i_\et(\Spa(R_j[\f{1}{p}], R_j), G) \qq \text{for every $i$.}
\ee
Since \cite{Hub96}*{3.2.9} applies to each $(R_j[\f{1}{p}], R_j)$ and gives a natural identification
\[
\tst H^i_\et(\Spa(R_j[\f{1}{p}], R_j), G) \cong H^i_\et(R_j^h[\f{1}{p}], G)
\]
(the definition of $\Spa(R_j[\f{1}{p}], R_j)_\et$ that we are using agrees with the one in \emph{op.~cit.},~see \cref{foot} and \cite{Hub96}*{2.2.8}), the combination of \eqref{Hi-alg-lim} and \eqref{Hi-an-lim} gives \eqref{iso-1}. The proof for \eqref{iso-3} is analogous: after expressing $R^\flat$ as a filtered direct limit of finite type $\bF_p$-subalgebras $R_j'$ that contain $\varpi^\flat$ and are integrally closed in $R_j'[\f{1}{\varpi^\flat}]$ and Henselizing each $R_j'$ with respect to the ideal $(\varpi^\flat) \subset R_j'$, one repeats the same arguments.
\epf

}



\section{Passage to perfect or perfectoid towers}

\ready{
We are ready to combine the results of the previous sections into a proof of the remaining cases of the purity conjecture for the Brauer group. We begin with auxiliary lemmas that build suitable~towers.

\blem \lab{mend-k}
For a complete, regular, local ring $(R, \fm)$, there is a filtered direct system of finite, flat $R$-algebras $R_i$ such that each $(R_i, \fm R_i)$ is a regular local ring and $(\varinjlim_i R_i, \fm(\varinjlim_i R_i))$ is a regular local ring with an algebraically closed residue field\uscolon if the residue field $k$ of $R$ is separably closed of characteristic $p > 0$, then the $R_i$ may be chosen to be of $p$-power rank over $R$.
\elem

\bpf
We set $p \ce \Char k$ and begin with the case when $R$ is of equicharacteristic. By the Cohen structure theorem \cite{Mat89}*{29.7}, then $R \simeq k\llb x_1, \dotsc, x_n\rrb$. We let $k_i$ range over the finite subextensions of $\ov{k}/k$ and set $R_i \ce k_i\llb x_1, \dotsc x_n\rrb$. The $\fm$-adic completion of $\varinjlim_i R_i$ is $\ov{k}\llb x_1, \dotsc x_n\rrb$, so $\varinjlim_i R_i$ is Noetherian (see \cite{SP}*{\href{http://stacks.math.columbia.edu/tag/033E}{033E}, \href{http://stacks.math.columbia.edu/tag/05UU}{05UU}, \href{http://stacks.math.columbia.edu/tag/00MK}{00MK}}), 
and hence also a regular local ring.

Now we turn to the case when $R$ is of mixed characteristic and $p \in \fm\setminus \fm^2$. Then, by \cite{Mat89}*{29.7} again, $R \simeq W\llb x_1, \dotsc, x_n \rrb$ for some complete discrete valuation ring $W$ with $p$ as a uniformizer. By \cite{Mat89}*{proof of 29.1}, there is an integral extension $W'/W$ of discrete valuation rings such that $W'$ has $p$ as a uniformizer and $\ov{k}$ as the residue field. Letting $W_i/W$ range over the finite discrete valuation ring subextensions of $W'/W$, we argue as in the equicharacteristic case that the local ring 
\[
\tst \p{\varinjlim_i (W_i\llb x_1, \dotsc, x_n\rrb), \fm(\varinjlim_i (W_i\llb x_1, \dotsc, x_n\rrb))}
\]
has $\wh{W'}\llb x_1, \dotsc, x_n\rrb$ as its completion and is regular.

In the remaining case when $R$ is of mixed characteristic and $p \in \fm^2$, by \cite{Mat89}*{29.3 and the proof of 29.8~(ii)}, 
there is a $W$ as above such that 
\[
R \simeq W\llb x_1, \dotsc, x_n\rrb/(p - f) \q \x{with} \q f \in (p, x_1, \dotsc, x_n)^2.
\]
Then, with the same $W'$ and $W_i$ as before, each 
\[
R_i \ce W_i\llb x_1, \dotsc, x_n\rrb/(p - f)
\]
is a finite, flat $R$-algebra that is a regular local ring. In addition, by the previous case, 
\[
\tst \varinjlim_i (W_i\llb x_1, \dotsc, x_n\rrb)
\]
is a regular local ring with a regular system of parameters $(p, x_1, \dotsc, x_n)$, so $\varinjlim_i R_i$ is a regular local ring with a regular system of parameters $(x_1, \dotsc, x_n)$.
\epf

The following variant of \cite{And18}*{3.4.5 (2)} supplies the perfectoid covers we will need.

\blem \lab{build-cover}
For a mixed characteristic $(0, p)$, complete, regular, local ring $(R, \fm)$ with a perfect residue field $k$, there is a tower $\{ R_m\}_{m \in \bZ_{\ge 0}}$ of finite, flat $R$-algebras $R_m$ of $p$-power rank over $R$ such that each $R_m$ is a regular, local ring and the $p$-adic completion $\wh{R}_\infty$ of $R_\infty \ce \varinjlim_m R_m$ is perfectoid. 
\elem

\bpf
We begin with the case when $p \in \fm \setminus \fm^2$, in which, by the proof of \Cref{mend-k} and the perfectness of $k$, we have $R \simeq W\llb x_1, \dotsc, x_n \rrb$ with $W \cong W(k)$. We set 
\[
R_m \ce (W[p^{1/p^m}])\llb x_1^{1/p^m}, \dotsc, x_n^{1/p^m}\rrb
\]
and use \Cref{p-oid-rec} to confirm that the resulting $\wh{R}_\infty$ is perfectoid.

In the remaining case when $p \in \fm^2$, we have $R \simeq W\llb x_1, \dotsc, x_n\rrb/(p - f)$ with $W \cong W(k)$ and some $f \in (p, x_1, \dotsc, x_n)^2$, and we set 
\[
R_m \ce W\llb x_1^{1/p^m}, \dotsc, x_n^{1/p^m}\rrb/(p - f).
\]
Due to \Cref{p-oid-rec}, to show that the $p$-adic completion of the local ring $(R_\infty, \fm_\infty)$ is perfectoid, we only need to argue that for some $u \in R^\times_\infty$ the element $up$ is a $p$-th power in $R_\infty$. For this, we follow the argument of \cite{Shi16}*{4.9} (alternatively, we could apply \cite{GR17}*{17.2.14~(ii)}): every element of $R_\infty/p$ is a $p$-th power and $p \in \fm_\infty^2$, so, by writing 
\[
\tst p = \sum_i ((s_i^p + pt_i)(s_i^{\prime p} + pt'_i)) \q \x{with} \q s_i, s_i' \in \fm_\infty \q \x{and}\q t_i, t_i' \in R_\infty,
\]
we find that $p = s^p + pt$ with $s, t \in \fm_\infty$ and may set $u \ce 1 - t$.
\epf

The key purity conclusion for the Brauer group is the following result.

\bthm \lab{win}
For a strictly Henselian, regular, local ring $(R, \fm)$ of dimension $\ge 2$, we have
\[
H^2_\et(U_R, \bG_m) = 0.
\]
\ethm

\bpf
We set $k \ce R/\fm$ and $p \ce \Char k$ and use \Cref{H2-inj} to assume that $R$ is complete. The case $\dim R = 2$ follows from \cite{Gro68c}*{6.1 b)} and the case $\dim R = 3$ then follows from \cite{Gab81}*{I, Thm.~2}, so we assume further that $\dim R \ge 4$. We then combine \Cref{cor-H2-inj} with a limit argument and \Cref{mend-k} to reduce to the case when $k = \ov{k}$ (to preserve completeness, we again use \Cref{H2-inj}). The absolute purity conjecture proved by Gabber, more precisely, \cite{Fuj02}*{2.1.1}, implies that for every prime $\ell \neq p$ we have $H^2_\et(U_R, \mu_\ell) = 0$, so also $H^2_\et(U_R, \bG_m)[\ell] = 0$. Therefore, since $H^2_\et(U_R, \bG_m)$ is torsion (see \Cref{reg-inj}), we will focus on the vanishing of $H^2_\et(U_R, \bG_m)[p]$.

We may assume that $p > 0$ and begin with the case when $R$ is an $\bF_p$-algebra, so that, as in the proof of \Cref{mend-k}, we have $R \simeq k\llb x_1, \dotsc, x_n \rrb$. Since $k = \ov{k}$, the $p^m$-Frobenius of $R$ is finite and flat for every $m > 0$, so, by combining \Cref{cor-H2-inj} with a limit argument, we reduce to proving that the perfection $V$ of $U_R$ satisfies $H^2_\et(V, \bG_m)[p] = 0$. This, in turn, is a special case of \Cref{perf-Gm}.

In the remaining case when $R$ is of mixed characteristic $(0, p)$, let $\{R_m\}$ be a tower supplied by \Cref{build-cover}. By \Cref{cor-H2-inj} and a limit argument, it suffices to show that $H^2_\et(U_{R_\infty}, \bG_m)[p] = 0$.  Each $R_m$ is regular, so, by \Cref{reg-inj},
\[
\tst H^2_\et(U_{R_\infty}, \bG_m) \hra H^2_\et(R_\infty[\f{1}{p}], \bG_m), \q \text{and, by \Cref{GR-inj},} \q  H^2_\et(R_\infty[\f{1}{p}], \bG_m) \hra H^2_\et(\wh{R}_\infty[\f{1}{p}], \bG_m).
\]
However, by \Cref{p-oid-Br}, the group $H^2_\et(\wh{R}_\infty[\f{1}{p}], \bG_m)$ has no nonzero $p$-torsion.
\epf

\bpp[Proof of \Cref{main-thm}] \lab{main-thm-pf}
We have a Henselian, regular, local ring $R$ of dimension $\ge 2$ whose residue field is of dimension $\le 1$ and an $R$-torus $T$, and we need to show that
\[
H^1_\et(U_R, T) = H^2_\et(U_R, T) = 0.
\]
By passing to the limit over all the finite, \'{e}tale, local $R$-algebras $R'$ and using \Cref{cor-H2-et}, we reduce to the case when $R$ is strictly Henselian. In this case, $T$ is split, $H^1_\et(U_R, \bG_m) = 0$ because every line bundle on $U_R$ extends to $R$, and $H^2_\et(U_R, \bG_m) = 0$ by \Cref{win}.
\QED
\epp

}



\section{Global conclusions}

The following standard arguments deduce  the global results stated in \S\ref{intro} from \Cref{win}. We begin with a mild generalization of \Cref{purity-conj} that simultaneously reproves \Cref{main-thm}:

\bthm\lab{purity-conj-pf}
For a scheme $X$, an $X$-torus $T$, and a closed subscheme $Z \subset X$ such that for every $z \in Z$ the local ring $\cO_{X,\, z}$ of $X$ at $z$ is regular of dimension $\ge 2$ and the inclusion $X \setminus Z \hra X$ is quasi-compact, we have 
\[
H^q_\et(X, T) \isomto H^q_\et(X - Z, T) \q \text{for} \q q \le 2 \qq \text{and} \qq H^3_\et(X, T) \hra H^3_\et(X - Z, T).
\]
\ethm

\bpf
By \cite{SGA4II}*{V, 6.5}, for each $X$-\'{e}tale $X'$ and the preimage $Z' \subset X'$ of $Z$, we have the exact sequence
\be \lab{coho-sup-seq}
\dotsc \ra H^2_{Z'}(X', T) \ra H^2_\et(X', T) \ra H^2_\et(X' - Z', T) \ra H^3_{Z'}(X', T) \ra \dotsc,
\ee
so it suffices to show that $H^q_Z(X, T) = 0$ for $q \le 3$. Thus, letting $\cH^q_Z(-, T)$ denote the \'{e}tale sheafification of the presheaf $X' \mapsto H^q_{Z'}(X', T)$, the local-to-global $E_2$ spectral sequence 
\[
H^p_\et(X, \cH^q_Z(X, T)) \Rightarrow H^{p + q}_Z(X, T)
\]
of \cite{SGA4II}*{V, 6.4} reduces us to showing that $\cH^q_Z(-, T) = 0$ for $q \le 3$. This is a local question, so, since tori are \'{e}tale locally trivial (see \cite{SGA3II}*{X, 4.5}), we may assume that $T = \bG_m$ and use a limit argument \cite{SGA4II}*{VII, 5.9} and the sequences \eqref{coho-sup-seq} to identify the stalk of $\cH^q_Z(-, \bG_m)$ at a geometric point $\ov{x}$ of $X$ with $H^q_{Z \times_{X} \Spec(\cO_{X,\, \ov{x}}^\sh)}(\cO_{X,\, \ov{x}}^\sh, \bG_m)$. Thus, we have reduced the desired vanishing of $H^{q}_Z(X, \bG_m)$ for $q \le 3$ to the case when $X$ is the spectrum of a strictly Henselian regular local ring of dimension $\ge 2$ and $Z\neq \emptyset$, so that $X$ is Noetherian and finite dimensional.

Therefore, we now assume that $X$ is Noetherian and finite dimensional and use the coniveau spectral sequence of \cite{Gro68c}*{\S10.1} (see also \cite{ILO14}*{XVIII$_{\x{{\upshape A}}}$, 2.2.1}):
\[
E_1^{pq} = \bigoplus_{z \in Z\x{ with } \dim(\cO_{X,\, z}) = p} H^{p + q}_{\{z\}}(\cO_{X,\, z}, \bG_m) \Ra H^{p + q}_Z(X, \bG_m),
\]
which reduces us further to the case when $Z$ is the closed point of a regular local  $X$. By the previous paragraph, we may then assume that $X$ is even strictly local (and regular, of dimension $\ge 2$). In this case, the vanishing $H^q_Z(X, \bG_m) = 0$ for $q \le 3$ follows from the bijectivity of the map $H^0(X, \bG_m) \isomto H^0_\et(X \setminus Z, \bG_m)$, the vanishing of $H^i_\et(X, \bG_m)$ for $i > 0$, and the vanishing of $H^i_\et(X \setminus Z, \bG_m)$ for $i = 1, 2$ supplied by the extendibility of line bundles and \Cref{win}. 
\epf

We turn to a more general version of \Cref{AG-cor}:

\bthm\lab{AG-cor-pf}
For a Noetherian, integral, regular scheme $X$ with function field $K$ and an $X$-torus $T$ such that the pullback map 
\be \lab{inj-ass}
H^2_\et(U, T) \ra H^2_\et(K, T) \qq \text{is injective for every nonempty open} \qq U \subset X
\ee
\up{for instance, $T$ could be $\bG_m$, see Lemma \uref{reg-inj}\uscolon for further examples, see Remark \uref{flasque}}, we have
\[
\tst H^2_\et(X, T) = \bigcap_{x \in X\,\text{of height}\, 1} H^2_\et(\cO_{X,\,x}, T) \qq \text{inside} \qq H^2_\et(K, T).
\]
\ethm

\bpf
Let $\gA$ be an element in the intersection. If $U, V \subset X$ are nonempty open subschemes such that $\gA$ extends to an element of both $H^2_\et(U, T)$ and $H^2_\et(V, T)$, then $\gA$ extends to an element of $H^2_\et(U \cup V, T)$: indeed, due to the injectivity assumption, this follows from the Mayer--Vietoris sequence
\be \lab{MV-seq}
\tst \dotsc \ra H^2_\et(U \cup V, T) \ra H^2_\et(U, T) \oplus H^2_\et(V, T) \ra H^2_\et(U \cap V, T) \ra H^3_\et(U \cup V, T) \ra \dotsc
\ee
that results from the \v{C}ech-to-derived spectral sequence $\check{H}^p(\{U, V\}, H^q_\et(-, T)) \Rightarrow H^{p + q}_\et(U \cup V, T)$ for the cover $\{U, V\}$ of $U \cup V$. 
Thus, $\gA$ extends to an element of $H^2_\et(U, T)$ for some nonempty open $U \subset X$ that covers the height $1$ points of $X$. Then, by \Cref{purity-conj-pf}, it also extends to $H^2_\et(X, T)$. 
\epf

\brems
\remi \lab{flasque}
By \cite{CTS87}*{2.2 (ii)}, the injectivity assumption \eqref{inj-ass} holds for any $T$ that is \emph{flasque} in the sense that there is a finite, \'{e}tale, Galois (so also connected, and hence nonempty) cover $X'/X$ that splits $T$ such that the $\Gal(X'/X)$-module $\Hom_{X'\text{-}\gp}(T_{X'}, \bG_{m,\, X'})$ has no nontrivial $\Gal(X'/X)$-module extensions by any module of the form $\bZ[\Gal(X'/X)/H]$ for a subgroup $H \le \Gal(X'/X)$ (see \cite{CTS87}*{0.5 and 1.2}). For example, any torus direct factor of $\Res_{X''/X}(\bG_m)$ for some finite \'{e}tale cover $X''/X$ is flasque.

On the other hand, the injectivity of $H^2_\et(X, T) \ra H^2_\et(K, T)$ fails for some $X$ and $T$. For example, consider the base change to $X \ce \bP^1_{\bR}$ of the torus $T$ defined by the exact sequence
\be \lab{SES-tors}
\qq 0 \ra T \ra \Res_{\bC/\bR}(\bG_m) \xra{\Norm} \bG_{m,\, \bR} \ra 0.
\ee
This sequence may be viewed as a $T$-torsor over $\bG_{m,\, \bR} \subset X$ whose $\bR$-fibral classes sweep out $H^1(\bR, T) \cong \bZ/2\bZ$. In particular, the torsor is not constant, that is, it is not a base change to $\bG_{m,\, \bR}$ of a $T$-torsor over $\bR$. In contrast, by \cite{CTS87}*{2.4}, every $T$-torsor over $\bA^1_\bR$ is constant. Consequently, due to the Mayer--Vietoris sequence \eqref{MV-seq}, the $T$-torsor defined by \eqref{SES-tors} gives rise to a nonzero, generically trivial class in $H^2_\et(\bP^1_\bR, T)$.

\remi
For a discrete valuation ring $\cO$ with the fraction field $K$ and the residue field $k$, by \cite{Gro68c}*{2.1}, for every prime $\ell \neq \Char k$, there is an exact residue sequence
\[
\qq 0 \ra H^2_\et(\cO, \bG_m)[\ell^\infty] \ra H^2_\et(K, \bG_m)[\ell^\infty] \ra H^1_\et(k, \bQ/\bZ)[\ell^\infty],
\]
and likewise for $\ell = \Char k$ if $k$ is perfect (but not in general when $\ell = \Char k$, see \cite{Poo17}*{6.8.2}). Therefore, \Cref{AG-cor-pf} implies that, as predicted in \cite{Poo17}*{6.8.4}, for a Noetherian, integral, regular scheme $X$ with the function field $K$, there is an exact sequence
\be \lab{res-seq}
\tst \qq 0 \ra H^2_\et(X, \bG_m) \ra H^2_\et(K, \bG_m) \ra \bigoplus_{x \in X\,\text{of height}\,1} H^1_\et(k(x), \bQ/\bZ)
\ee
granted that one excludes the $p$-primary parts for every prime $p$ for which some point $x \in X$ of height $1$ has an imperfect residue field $k(x)$ of characteristic $p$.\footnote{In particular, one excludes the $p$-primary parts whenever $X$ has a point $y \in X$ of residue characteristic $p$ and height $\ge 2$. Indeed, such a $y$ generalizes to some height $1$ point $x \in X$ of residue characteristic $p$ and, by \cite{EGAII}*{7.1.7}, the residue field $k(x)$ has a nontrivial discrete valuation, so is imperfect. Thus, the special cases of \Cref{AG-cor-pf} that had been known prior to this work suffice for deducing \eqref{res-seq}.}

\erems


\begin{appendix}

\section{Fields of dimension $\le 1$}

\ready{
The formulation of \Cref{main-thm} above involves the following well-known class of fields.

\bd[\cite{Ser02}*{II.\S3.1, Prop.~5 and I.\S3.1, Prop.~11}] \lab{dim-1-def}
A field $k$ is \emph{of dimension $\le 1$} if 
\[
H^i_\et(k, G) = 0 \q \text{for} \q i \ge 2 \q \text{and every commutative, finite, \'{e}tale $k$-group scheme $G$} 
\]
and if also, when $\Char k$ is positive, $H^2_\et(K, \bG_m) = 0$ for every finite, separable extension $K/k$.
\ed

In this short appendix, we record \Cref{dim-1-tori} in the form convenient for its use in the proof of \Cref{cor-H2-et} and give an equivalent definition of a field of dimension $\le 1$ in \Cref{main}, which, we believe, deserves to be known more widely.

\blem \lab{dim-1-tori}
For a field $k$ of dimension $\le 1$ and a $k$-torus $T$, we have
\be \lab{D1T-eq}
H^i_\et(k, T) = 0 \q \text{for every} \q i \ge 1.
\ee
\elem

\bpf
The strict cohomological dimension of $k$ is $\le 2$ (see \cite{Ser02}*{I.\S3.2, Prop.~13}), so the $i \ge 3$ case follows. Thus, so does the case $T = \bG_m$. Consequently, by choosing a finite Galois extension $K/k$ that splits $T$ and considering the norm map $\Res_{K/k} (T_K) \ra T$, at the cost of changing $T$, we may replace $H^i_\et$ by $H^{i + 1}_\et$ in \eqref{D1T-eq}, and then likewise by $H^{i + 2}_\et$. Thus, the settled $i \ge 3$ case suffices.
\epf

\bthm \lab{main}
A field $k$ is of dimension $\le 1$ if and only if 
\[
H^i_\fppf(k, G) = 0 \q \text{for} \q i \ge 2 \q \text{and every commutative, finite $k$-group scheme $G$.}
\]
\ethm

\bpf
We may assume that $p \ce \Char k$ is positive. The displayed condition implies that $k$ is of dimension $\le 1$: indeed, for $K/k$ finite, separable, each $H^2_\et(K, \bG_m) \cong H^2_\et(k, \Res_{K/k}(\bG_m))$ is torsion, and hence vanishes as soon as $H^2_\fppf(k, (\Res_{K/k}(\bG_m))[\ell])$ vanishes for every prime $\ell$ (including $\ell = p$).

For the converse, we assume that $k$ is of dimension $\le 1$ and, by decomposing and filtering $G$, that $G$ is killed by $p$, connected, and with $G^\vee$ that is either connected or \'{e}tale. If $G^\vee$ is also connected, then, by \cite{SGA3II}*{XVII, 4.2.1 ii) $\Leftrightarrow$ iv)}, the group $G$ is a successive extension of the copies of the Frobenius kernel $\gA_p$ of $\bG_a$. The vanishing of the coherent cohomology $H^i(k, \bG_a) = 0$ for $i \ge 1$ then gives the claim. If $G^\vee$ is \'{e}tale, then $G$ is the kernel of a map of $k$-tori and \Cref{dim-1-tori} suffices.
\epf

\bcor
A field $k$ is of dimension $\le 1$ if and only if 
\[
H^i_\fppf(k, G) = 0 \q \text{for} \q i \ge 2 \q \text{and every commutative, finite type $k$-group scheme $G$.}
\]
\ecor

\bpf
We may focus on the `only if.' In addition, by \cite{SGA3Inew}*{VII$_{\text{{\upshape A}}}$, 8.3} and \Cref{main}, we may assume that $G$ is smooth, and then also connected. Then $H^i_\fppf(k, G) \cong H^i_\et(k, G)$ is torsion for $i \ge 1$, so consideration of the $\ell$-torsion $G[\ell]$ settles the case when $\Char k = 0$ or $G$ is semiabelian. Thus, the ``anti-Chevalley theorem'' \cite{CGP15}*{A.3.9} reduces further to affine $G$. Grothendieck's theorem on maximal tori \cite{SGA3II}*{XIV, 1.1} then allows us to assume that $G$ is unipotent (see \cite{SGA3II}*{XVII, 4.1.1}). 
For unipotent $G$, the filtration of \cite{SGA3II}*{XVII, 3.5 ii)} suffices.
\epf

\brem 
For vanishing results for $H^1(k, G)$ with $k$ of dimension $\le 1$, see \cite{Ser02}*{III.\S2.3}.
\erem

}





\end{appendix}

\end{document}